\documentclass[11pt]{amsart}

\usepackage{amsmath}
\usepackage{amssymb}
\usepackage{graphicx}
\usepackage{amscd}
\usepackage{stmaryrd}
\usepackage[mathscr]{euscript}

\theoremstyle{definition}

\numberwithin{equation}{section}

\newcommand{\sB}{\mathscr B}
\newcommand{\sC}{\mathscr C}

\newcommand{\sE}{\mathscr E}

\newcommand{\sH}{\mathscr H}
\newcommand{\sJ}{\mathscr J}
\newcommand{\sK}{\mathscr K}
\newcommand{\sM}{\mathscr M}

\newcommand{\sR}{\mathscr R}
\newcommand{\sS}{\mathscr S}

\newcommand{\bR}{\mathbb R}

\title[Spectral Measures on Fractals]{The formula for the quasicentral modulus in the case of spectral measures on fractals}

\author[D.-V. Voiculescu]{Dan-Virgil Voiculescu${}^*$}

\address{Department of Mathematics \\ University of California at Berkeley \\ Berkeley, CA\ \ 94720-3840}
\email{{\tt dvv@math.berkeley.edu}}

\begin{document}

\begin{abstract}
We prove a general ampliation homogeneity result for the quasicentral modulus of an $n$-tuple of operators with respect to the $(p,1)$ Lorentz normed ideal. We use this to prove a formula involving Hausdorff measure for the quasicentral modulus of $n$-tuples of commuting Hermitian operators the spectrum of which is contained in certain Cantor-like self-similar fractals.
\end{abstract}

\maketitle


\section{Introduction}
\label{sec1}

The quasicentral modulus $k_{\sJ}(\tau)$ is a number associated with an $n$-tuple $\tau$ of Hermitian operators relative to a normed ideal $(\sJ,|\ |_{\sJ})$ of compact operators. It underlies many questions on normed ideal perturbations of $n$-tuples of operators (see the recent survey \cite{11}), and it also had applications in non-commutative geometry in work on the spectral characterization of manifolds \cite{1}.

We proved in \cite{7} that in the case of $\tau$ an $n$-tuple of commuting Hermitian operators and if the normed ideal is the $(n,1)$-Lorentz ideal, which we denote by $\sC_n^-$, the corresponding quasicentral modulus $k^-_n(\tau)$ has the property that $(k_n^-(\tau))^{1/n}$ is proportional to the integral w.r.t.\ $n$-dimensional Lebesgue measure of the multiplicity function of $\tau$.

Here we prove a similar result in fractional dimension. More precisely instead of a cube in $\bR^n$ which contains the spectrum of $\tau$ we assume there is a fixed self-similar fractal in $\bR^n$ of Hausdorff dimension $p > 1$ containing the spectrum $\sigma(\tau)$. The analogous formula we prove has the exponent $n$ replaced by $p$ and the integral of the multiplicity function is with respect to $p$-Hausdorff measure. For technical reasons the class of fractals is rather restricted, only certain totally disconnected sets, that is Cantor-like fractals are considered. One should certainly expect this can be extended to a larger class of fractals, the present paper being only a first step.

\bigskip
\noindent
\underline{\qquad\qquad} 

2020 {\em Mathematics Subject Classification}. Primary: 47A55; Secondary: 28A78, 28A80, 47B10.

{\em Key words and phrases}. Lorentz (p, 1) normed ideal, p-Hausdorff measure, quasicentral modulus.

\vfill
\noindent
${}^*$Research supported in part by NSF-Grant DMS-1665534.

\newpage
What made the extension of the formula to fractional dimension possible is a completely general ampliation homogeneity result $k^-_p(\tau \otimes I_m) = m^{1/p}k^-_p(\tau)$. Such a result was previously known only for $p = \infty$ (\cite{9}) and for $p = 1$ (\cite{7}). Note that in \cite{7} we obtained easily an ampliation homogeneity result for $k_p(\tau)$, that is with $\sC^-_p$ replaced by the Schatten--v.~Neumann class $\sC_p$. For $p = 1$ we have $\sC_1 = \sC^-_1$, but for $p > 1$ the result for $k_p(\tau)$ turned out to be trivial after we showed in \cite{8} that in this case $k_p(\tau) \in \{0,\infty\}$. The interesting quantity which replaces $k_p(\tau)$ is $k_p^-(\tau)$.

The paper has six sections including this introduction. Section~2 contains preliminaries about the quasicentral modulus. Section~3 is devoted to the ampliation homogeneity theorem. In section~4 we collected preliminaries concerning the class of Cantor-type fractals we consider. The formula for $k_p^-(\tau)$ in the fractal setting is obtained in section~5. Section~6 deals with concluding remarks.

\section{Operator preliminaries}
\label{sec2}

By $\sH$ we denote a separable complex Hilbert space of infinite dimension and by $\sB(\sH)$, $\sK(\sH)$, $\sR(\sH)$ the bounded operators, the compact operators and the finite rank operators. When no confusion will arise we will simply write $\sK,\sR$ and we will denote by $\sR_1^+(\sH)$ or $\sR_1^+$ the finite rank positive contractions $0 \le A \le I$ on $\sH$. The $(p,1)$ Lorentz normed ideal of compact operators will be denoted by $(\sC^-_p,|\ |_p^-)$. We recall that the norm is 
\[
|T|^-_p = \sum_{j \in {\mathbb N}} s_j\ j^{-1+1/p}
\]
where $s_1 \ge s_2 \ge \dots$ are the eigenvalues of $|T| = (T^*T)^{1/2}$ in decreasing order. If $(\sC_p,|\ |_p)$ is the Schatten--von~Neumann $p$-class, then $\sC_1 = \sC_1^-$. More on normed ideals can be found in \cite{5} and \cite{6}.

We shall also use the following notation for operations on $n$-tuples of operators, in line with \cite{7}. If $\tau = (T_i)_{1 \le i \le n} \in (\sB(\sH))^n$ and $X,Y \in \sB(\sH)$ then we use
\[
\begin{aligned}
X\tau Y &= (XT_iY)_{1 \le i \le n} \\
[X,\tau] &= ([X,T_i])_{1 \le i \le n} \\
\tau^* &= (T_i^*)_{1 \le i \le n}.
\end{aligned}
\]
If also $\sigma = (S_i)_{1 \le i \le n} \in (\sB(\sH))^n$ then we write
\[
\begin{aligned}
\sigma+\tau &= (S_i+T_i)_{1 \le i \le n} \\
\sigma\oplus\tau &= (S_i \oplus T_i)_{1\le i \le n} \\
\tau\otimes I_m &= (T_i \otimes I_m)_{1 \le i \le n}
\end{aligned}
\]
where $I_m$ is the identity operator on ${\mathbb C}^m$. When we identify $\sH \otimes {\mathbb C}^m$ and $\sH \oplus \dots \oplus \sH$ we also have
\[
\tau \otimes I_m \cong \underset{\text{$m$-times}}{\underbrace{\tau \oplus \dots \oplus \tau}}\,.
\]
Further, we consider norms
\[
\begin{aligned}
\|\tau\| &= \max_{1 \le i \le n} \|T_i\| \\
|\tau|_{\sJ} &= \max_{1 \le i \le n} |T_i|_{\sJ}.
\end{aligned}
\]

The quasicentral modulus of an $n$-tuple $\tau = (T_i)_{1 \le i \le n}$ with respect to a normed ideal $(\sJ,|\ |_{\sJ})$ (see \cite{7}, \cite{9}) is the number
\[
k_{\sJ}(\tau) = \liminf_{A \in \sR_1^+} |[\tau,A]|_{\sJ}
\]
where the $\liminf$ is w.r.t.\ the natural order on $\sR_1^+$. This definition is also equivalent to $k_{\sJ}(\tau)$ being the
\[
\begin{aligned}
\inf\{\alpha \in [0,\infty]| \alpha &= \lim_{k \to \infty} |[A_k,\tau]|_{\sJ},\ A_k \uparrow I \\
&A_k \in \sR_1^+\}
\end{aligned}
\]
or the same with $w - \lim_{k \to \infty} A_k = I$ instead of $A_k \uparrow I$. If $\sJ = \sC_p^-$ we denote $k_{\sJ}(\tau)$ by $k_p^-(\tau)$.

We should also record as the next proposition the results in \cite{7} Prop.~$1.4$ and Prop.~$1.6$.

\medskip
\noindent
{\bf Proposition 2.1.}
{\em If $\tau^{(j)} \in \sB(\sH)^n$, $j \in {\mathbb N}$ then we have
\[
\max_{j=1,2} k_{\sJ}(\tau^{(j)}) \le k_{\sJ}(\tau^{(1)} \oplus \tau^{(2)}) \le k_{\sJ}(\tau^{(1)}) + k_{\sJ}(\tau^{(2)})
\]
and
\[
k_{\sJ}\Big( \underset{j \in {\mathbb N}}{\bigoplus} \tau^{(j)}\Big) = \lim_{m_i \to \infty} k_{\sJ}\Big( \underset{1 \le j \le m}{\bigoplus} \tau^{(j)}\Big).
\]
If $\lambda^{(j)} \in {\mathbb C}^n$ and $\lambda^{(j)} \otimes I_{\sH} \in \sB(\sH)^n$ then we have
\[
\begin{aligned}
k_{\sJ}(\tau^{(1)} \oplus \dots \oplus \tau^{(m)}) &= \\
&= k_{\sJ}\left((\tau^{(1)}-\lambda^{(1)} \otimes I_{\sH}) \oplus \dots \oplus (\tau^{(m)} - \lambda^{(m)} \otimes I_{\sH})\right).
\end{aligned}
\]
}

\medskip
Finally, if $\tau$ is a $n$-tuple of commuting Hermitian operators we denote by $\sigma(\tau) \subset {\mathbb R}^n$ the joint spectrum and by $E(\tau;\omega)$ the spectral projection of $\tau$ for the Borel set $\omega \subset {\mathbb R}^n$.

\section{Ampliation homogeneity}
\label{sec3}

\medskip
\noindent
{\bf Theorem 3.1.} 
{\em If $\tau$ is a $n$-tuple of bounded operators and $1 \le p \le \infty$ then
\[
k_p^-(\tau \otimes I_m) = m^{1/p} k_p^-(\tau).
\]
}

\medskip
The cases $p = 1$ and $p = \infty$ have already been proved (\cite{7} Prop.~1.5 and \cite{9} Prop.~3.9).

We begin the proof with a couple of lemmas.

\medskip
\noindent
{\bf Lemma 3.1.} 
{\em Let $X_j \in {\mathcal C}^-_p$, $j \in {\mathbb N}$, $p \in [1,\infty]$ be so that $|X_j|^-_p \le C$. If $\lim_{j \to \infty} \|X_j\| = 0$ then we have
\[
\lim_{j \to \infty} (|X_j \otimes I_m|^-_p - m^{1/p}|X_j|_p^-) = 0.
\]
}

\begin{proof}
Let $s_1^{(j)} \ge s_2^{(j)} \ge \dots$ be the eigenvalues of $(X_j^*X_j)^{1/2}$. Then
\[
\begin{aligned}
|X_j \otimes I_m|^-_p &= \sum_{k \in {\mathbb N}} s_k^{(j)}((m(k-1)+1)^{-1+1/p} + \dots + (mk)^{-1+1/p}) \\
&\ge m^{1/p} \sum_{k \in {\mathbb N}} s_k^{(j)}k^{-1+1/p} = m^{1/p}|X_j|^-_p.
\end{aligned}
\]
On the other hand, given $\epsilon > 0$ there is $N$ so that
\[
k \ge N \Rightarrow (m(k-1)+1)^{-1+1/p} + \dots + (mk)^{-1+1/p} \le (1+\epsilon)m^{1/p}k^{-1+1/p}.
\]
This gives
\[
\begin{aligned}
|X_j \otimes I_m|_p^- &\le (1+\epsilon)m^{1/p} \sum_{k \ge N} s_k^{(j)}k^{-1+1/p} + N\|X_j\| \\
&\le (1+\epsilon)m^{1/p}|X_j|^-_p + N\|X_j\|.
\end{aligned}
\]
Thus we have
\[
0 \le |X_j \otimes I_m|^-_p - m^{1/p}|X_j|^-_p \le \epsilon m^{1/p}|X_j|_p^- + N\|X_j\|.
\]
Since $\epsilon > 0$ is arbitrary and $\|X_j\| \to 0$ we get the desired result when $j \to \infty$.
\end{proof} 

\medskip
\noindent
{\bf Corollary 3.1.} 
{\em Let $X_j = (X_{ji})_{1 \le i \le n}$ be $n$-tuples of operators so that $|X_j|_p^- \le C$ where $p \in [1,\infty]$ and $\lim_{j \to \infty} \|X_j\| = 0$. Then we have 
\[
\lim_{j \in \infty}(m^{1/p} |X_j|_p^- - |X_j \otimes I_m|^-_p) = 0.
\]
}

\medskip
\noindent
{\bf Lemma 3.2.} 
{\em If $\tau \in (\sB(\sH))^n$ and $(\sJ,|\ |_{\sJ})$ is a normed ideal so that $k_{\sJ}(\tau) < \infty$, then there are $B_j \in \sR^+_1$ so that $B_j \uparrow I$ and
\[
\begin{aligned}
\lim_{j \to \infty} |[\tau,B_j]|_{\sJ} = k_{\sJ}(\tau) \\
\lim_{j \to \infty} \|[\tau,B_j]\| = 0.
\end{aligned}
\]
}

\begin{proof}
It suffices to show that given $\epsilon > 0$ and $P$ a finite rank Hermitian projector we can find $B \in \sR_1^+$ so that $B \ge P$ and
\[
\begin{aligned}
|[B,\tau]|_{\sJ} &\le k_{\sJ}(\tau) + \epsilon \\
\|[B,\tau]\| &\le \epsilon.
\end{aligned}
\]
Such that $B$ can be constructed as follows. We find recursively $P = P_1 \le P_2 \le P_3 \le \dots$ finite rank Hermitian projectors and $A_j \in \sR_1^+$ so that
\[
\begin{aligned}
&A_j \ge P_j,\ |[A_j,\tau]|_{\sJ} \le k_{\sJ}(\tau) + \epsilon \\
&\tau P_j = P_{j+1}\tau P_j,\ \tau^*P_j = P_{j+1}\tau^*P_j \\
&P_{j+1} \ge A_j.
\end{aligned}
\]
If we put $Q_j = P_{j+1} - P_j$ if $j \ge 1$ and $Q_0 = P_1 = P$ we have
\[
Q_r\tau Q_s \ne 0 \Rightarrow |r-s| \le 1
\]
and
\[
A_j = (Q_0 + \dots + Q_{j-1}) + Q_jA_jQ_j.
\]
This gives
\[
Q_r[\tau,A_j]Q_s \ne 0
\]
$\Rightarrow |r-s| \le 1$ and $j-1 \le r,s \le j+1$. It follows that if $|k-j| \ge 4$ $([\tau,A_j])^*[\tau,A_k] = 0$ so that
\[
\|[\tau,A_4 + A_8 + \dots + A_{4N}]\| \le 2\|\tau\|.
\]
Thus if $B_N = N^{-1}(A_4 + A_8 + \dots + A_{4N})$ we have $B_N \ge P$, $B_N \in \sR_1^+$, $|[B_N,\tau]|_{\sJ} < k_{\sJ}(\tau) + \epsilon$ and
\[
\|[B_N,\tau]\| \le 2N^{-1}\|\tau\|.
\]
Thus if $2N^{-1}\|\tau\| < \epsilon$ we may take $B = B_N$.
\end{proof}

\medskip
\noindent
{\bf Lemma 3.3.} 
{\em If $m \in {\mathbb N}$ and $(\sJ,|\ |_{\sJ})$ is a normed ideal and $\tau \in (B(\sH))^n$ is so that $k_{\sJ}(\tau) < \infty$, then there are $A_j \in \sR^+_1$, $A_j \uparrow I$ so that
\[
k_{\sJ}(\tau \otimes I_m) = \lim_{j \to \infty} |[\tau \otimes I_m,A_j \otimes I_m]|_{\sJ}
\]
and $\lim_{j \to \infty} \|[\tau,A_j]\| = 0$.
}

\begin{proof}
Let $G$ be the group $\sS_m \rtimes {\mathbb Z}_2^m$ of permutation matrices with $\pm 1$ entries and $g \to U_g$ its representation on $\sH \oplus \dots \oplus \sH$ which is $\otimes I_{\sH}$ the representation on ${\mathbb C}^m$. Then the commutant $\{U_g: g \in G\}'$ is $B(\sH) \otimes I_m$ and the map $\Phi: B(\sH^m) \to \sB(\sH) \otimes I_m$ given by $\Phi(X) = |G|^{-1} \sum_{g \in G} U_gXU_g^*$ is the projection of norm one which preserves the trace. If $B \in R_1^+(\sH^m)$ we have
\[
|[B,\tau \otimes I_m]|_{\sJ} = |U_g[B,\tau\otimes I_m]U_g^*|_{\sJ} = |[U_gBU_g^*,\tau \otimes I_m]|_{\sJ},
\]
which gives by taking the mean over $G$
\[
|[\Phi(B),\tau \otimes I_m]|_{\sJ} \le |[B,\tau \otimes I_m]|_{\sJ}
\]
and clearly $B_j \uparrow I \otimes I_m$ implies $\Phi(B_j) \uparrow I \otimes I_m$. Thus if $A_j \otimes I_m = \Phi(B_j)$ and $B_j \uparrow I \otimes I_m$ are so that
\[
\lim_{j \to \infty} |[B_j,\tau \otimes I_m]|_{\sJ} = k_{\sJ}(\tau \otimes I_m)
\]
and
\[
\lim_{j \to \infty} |[B_j,\tau \otimes I_m]\| = 0
\]
then we have
\[
\lim_{j \to \infty} \sup|[A_j \otimes I_m,\tau \otimes I_m]|_{\sJ} \le k_{\sJ}(\tau \otimes I_m)
\]
and
\[
\lim_{j \to \infty} \|[A_j \otimes I_,\tau \otimes I_m]\| = 0.
\]
Since on the other hand
\[
\liminf_{j \to \infty} |[A_j \otimes I_m,\tau \otimes I_m]|_{\sJ} \ge k_{\sJ}(\tau \otimes I_m)
\]
we conclude that
\[
\lim_{j \to \infty} |[A_j \otimes I_m,\tau \otimes I_m]|_{\sJ} = k_{\sJ}(\tau).
\]
\end{proof}

\medskip
\noindent
{\bf Proof of Theorem~3.1.} Using Lemma~$3.3$ we can find $A_j \in \sR_1^+$, $A_j \uparrow I$ when $j \to \infty$ so that
\[
\lim_{j \to \infty} |[A_j,\tau] \otimes I_m]|^-_p = k^-_p(\tau \otimes I_m)
\]
and
\[
\lim_{j \to \infty} \|[A_j,\tau]\| = 0.
\]
Using Corollary~$3.1$ we infer that
\[
\lim_{j \to \infty} m^{1/p}|[A_j;\tau]|^-_p = k_p^-(\tau \otimes I_m)
\]
which implies that
\[
m^{1/p}k^-_p(\tau) \le k^-_p(\tau \otimes I_m).
\]
On the other hand, Lemma~$3.2$ shows that there are $A_j \uparrow I$, $A_j \in \sR_1^+$ so that
\[
\lim_{j \to \infty} |[\tau,A_j]|^-_p = k^-_p(\tau)
\]
and
\[
\lim_{j \to \infty} \|[\tau,A_j]\| = 0.
\]
Then by Corollary~$3.1$ we get that
\[
\begin{aligned}
m^{1/p}k^-_p(\tau) &= \lim_{j \to \infty} |[\tau,A_j] \otimes I_m|^-_p \\
&= \lim_{j \to \infty} |[\tau \otimes I_m,A_j \otimes I_m]|^-_p \\
&\ge k^-_p(\tau \otimes I_m)
\end{aligned}
\]
which concludes the proof.\qed

\section{Fractal preliminaries}
\label{sec4}

To keep things simple the fractal context will be certain totally disconnected Cantor-like self-similar sets and on which a certain Hausdorff measure is a Radon measure on Borel sets.

We consider a non-empty compact set $K \subset \bR^n$ and a $N$-tuple of maps
\[
F_i(x) = \lambda (x - b(i)) + b(i),\ 1 \le i \le N
\]
where $0 < \lambda < 1$ and $b(i) \in \bR^n$, so that
\[
K = \bigcup_{1 \le i \le N} F_iK
\]
and we assume that
\[
i_1 \ne i_2 \Rightarrow F_{i_1}K \cap F_{i_2}K = \emptyset.
\]
Note that the open set condition for the contractions $F_i$ (see \cite{4}, page~121) in this case can be satisfied with an open neighborhood of $K$. The Hausdorff and box-dimension of $K$ are equal to
\[
p = \frac {\log N}{\log(1/\lambda)}
\]
by \cite{4} Thm.~$8.6$, and the $p$-Hausdorff measure of $K$ is finite and non-zero. This Hausdorff measure is often referred to as the Hutchinson measure of $K$. Note also the uniqueness of $K$ given the maps $F_i$, $1 \le i \le N$ (\cite{4} Thm.~$8.3$).

If $w \in \{1,\dots,N\}^m$ we put $|w| = m$ and define $F_w = F_{w_1} \circ \dots \circ F_{w_m}$ and $K_w = F_wK$. In particular if $|w| = |w'|$ then $K_w$ and $K_{w_1}$ are congruent and in particular have the same $p$-Hausdorff measure and diameter. Moreover
\[
K = \bigcup_{|w|=L} K_w
\]
and for any $L \in {\mathbb N}$ the union is disjoint. On $K$ there is a unique Radon measure $\mu$ so that
\[
\mu(K_w) = N^{-|w|} = \lambda^{|w|p}
\]
and on Borel sets $\mu = cH_p$ where $H_p$ denotes the $p$-Hausdorff measure and $c$ is a constant.  We use the $l^{\infty}$-norm $|(x_1,\dots,x_n)| = \max_{1 \le i \le n} |x_i|$ on $\bR^n$. Note that
\[
\mbox{diam}(K_w) = c \cdot \lambda^{|w|}
\]
for some constant $c$.

We shall also assume that the Hausdorff dimension $p \ge 1$.

If $\tau$ is a $n$-tuple of commuting Hermitian operators on $\sH$ with spectrum $\sigma(\tau) \subset K$, like with Lebesgue measure here on $K$ with $\mu$ that is with $H_p$ the Hilbert space splits
\[
\sH = \sH_{psing} \oplus \sH_{pac}
\]
where $\sH_{psing}$, $\sH_{pac}$ are reducing subspaces for $\tau$ and consist of vectors $\xi$ so that
\[
\langle E(\tau;\cdot)\xi,\xi\rangle
\]
is singular w.r.t.\ $\mu$ and respectively absolutely continuous w.r.t.\ $\mu$, that is w.r.t.\ $H_p$. (Use for instance $1.6.3$ in \cite{3}.) For more on Hausdorff-measure and on fractals see \cite{3} and \cite{4}.

\section{$k^-_p$ in the fractal setting}
\label{sec5}

In this section we study $k^-_p(\tau)$ where $\tau$ is a $n$-tuple of commuting Hermitian operators with $\sigma(\tau) \subset K$, in the context of Section~4. We assume $p \ge 1$ (and for certain results we will require $p > 1$).

\medskip
\noindent
{\bf Lemma 5.1.} 
{\em Assume $\tau$ is a $n$-tuple of commuting Hermitian operators with $\sigma(\tau) \subset K$ and with a cyclic vector $\xi$. Then for some constant $C$ depending only on $K$, we have
\[
k^-_p(\tau) \le C(H_p(\sigma(\tau)))^{1/p}.
\]
}

\begin{proof}
Let $\Omega(L) = \{w \mid |w| = L,\ K_w \cap \sigma(\tau) \ne \emptyset\}$ and $G(L) = \bigcup_{w \in \Omega(L)} K_w$. Then $G(L)$ is open in $K$ and there is $L_0$ so that for a given $\epsilon > 0$ we have $L \ge L_0 \Rightarrow H_p(G(L)) \le H_p(\sigma(\tau)) + \epsilon$. Let further $E_w = E(\tau;K_w)$ and observe that since $G(L) \supset \sigma(\tau)$ we have
\[
\sum_{w \in \Omega(L)} E_w = I.
\]
We also have
\[
H_p(G_L) = |\Omega(L)|H_p(K)\lambda^{L_p}.
\]
Let $P_L$ be the orthogonal projection onto $\sum_{w \in \Omega(L)} {\mathbb C}E_w\xi$. and $P_w$ the orthogonal projection onto ${\mathbb C}E_w\xi$. Remark that the non-zero $E_w\xi$ are an orthogonal basis on $P_L\sH$ and rank $P_L \le |\Omega_L|$. We have $[P_L,E_w] = 0$ if $|w| = L$ and $L_1 \ge L_2 \Rightarrow P_{L_1} \ge P_{L_2}$ since each  $E_w$ with $|w| = L_2$ is the sum of $E_{w'}$ with $|w'| = L_1$. Remark also that
\[
\begin{aligned}
\|[P_L,\tau]\| &= \max_{|w|=L} \|[P_L,\tau]E_w\| \\
&= \max_{|w| = L}\|[P_w,E_w\tau]\| \\
&\le \max_{|w|=L} 2\mbox{diam}(\sigma(E_w\tau \mid E_w\sH)) \\
&\le \max_{|w|=L} 2\mbox{diam}(K_2) \le 2\lambda^L\mbox{diam}(K).
\end{aligned}
\]
If $P_L \uparrow P$ then $P\xi = \xi$ and $[P,\tau] = 0$ since $\|[P_L,\tau]\| \to 0$. The vector $\xi$ being cyclic this gives $P = I$, so $P_L \uparrow I$ as $L \to \infty$. Hence we get
\[
\begin{aligned}
|[P_L,\tau]|^-_p &\le p(\mbox{rank}[P_L,\tau])^{1/p} \cdot \|[P_L,\tau]\| \\
&\le p(2|\Omega(L)|)^{1/p} \cdot 2\lambda^L\mbox{diam}(K) \\
&= p \cdot 2^{1+1/p} \cdot \left( \frac {H_p(G_L)}{H_p(K)} \lambda^{-Lp}\right)^{1/p} \cdot \lambda^L\mbox{diam } K \\
&\le p \cdot 2^{1+1/p}(H_p(K))^{-1/p} \cdot \mbox{diam } K(H_p(\sigma(\tau))+\epsilon)^{1/p}.
\end{aligned}
\]
Since $\epsilon > 0$ we get
\[
|[P_L,\tau]|^-_p \le c(H_p(\sigma(\tau)))^{1/p}
\]
where $c = p \cdot 2^{1+1/p}\cdot\mbox{diam } K(H_p(K))^{-1/p}$. Letting $L \to \infty$ we have
\[
k_p^-(\sigma(\tau)) \le c(H_p(\sigma(\tau)))^{1/p}.
\]
\end{proof}

\medskip
\noindent
{\bf Lemma 5.2.}
{\em Assume $\sigma(\tau) \subset K$ and assume that the spectral measure $E(\tau;\cdot)$ is singular w.r.t.\ $H_p$. Then we have
\[
k^-_p(\tau) = 0.
\]
}

\begin{proof}
Since $\tau$ is the orthogonal sum of $n$-tuples with cyclic vector, it suffices to prove the lemma when $\tau$ has a cyclic unit vector $\xi$. The absolute continuity class on $K$ of $E(\tau;\cdot)$ is then the same as the absolute continuity class of the scalar measure $\nu =\langle E(\tau;\cdot)\xi,\xi\rangle$.

Given $\epsilon > 0$ we can find a compact set $C_m$ which is a finite union of $K_w$ so that $H_p(C_m) \le \epsilon$ and $\nu(K\backslash C_m) < 2^{-m}$. Since $\sigma(\tau \mid E(\tau;C_m)\sH) \subset C_m$ the preceding lemma gives that
\[
\begin{aligned}
&k^-_p(\tau \mid E(\tau;C_m)\sH) \le \\
&\le c \cdot H_p(C_m) \le c \cdot \epsilon.
\end{aligned}
\]
On the other hand $\nu(K\backslash C_m) \le 2^{-m}$ gives
\[
\|\xi - E(\tau;C_m)\xi\|^2 \le 2^{-m}.
\]
Since $(\tau)''$ is a maximal abelian von~Neumann algebra in ${\mathcal B}({\mathcal H})$, $\xi$ being cyclic is also separating and we infer $E(\tau;C_m)\ {\overset{w}{\rightarrow}}\ I$ and hence 
$E(\tau;C_m)\ {\overset{s}{\rightarrow}}\ I$ as $m \to \infty$. This gives that we can find $A_m \in \sR_1^+$, $A_m \le E(\tau;C_m)$, $|[A_m,\tau]|^-_p \le k_p^-(\tau \mid E(\tau;C_m))$ and $A_m \uparrow I$. It follows that
\[
k^-_p(\tau) \le c \cdot \epsilon
\]
and $\epsilon$ being arbitrary $k^-_p(\tau) = 0$.
\end{proof}

Note that on $K$ the $p$-Hausdorff measure satisfies an Ahlfors regularity condition
\[
C^{-1}r^p \le H_p(B(x,r)) \le Cr^p
\]
if $r \le 1$ for some $C > 0$. The right half of this $H_p(B(x,r)) \le Cr^p$ if $r \le 1$ is the sub-regularity condition where $p > 1$ required in \cite{2} Cor.~$4.7$ to show that $k^-_p(\tau_K) > 0$ where $\tau_K$ is the $n$-tuple of multiplication operators by the coordinate functions in $L^2(K,H_p \mid K)$. Thus we have

\medskip
\noindent
{\bf Lemma 5.3.} ([2]) Assume $p > 1$ then $k^-_p(\tau_K) > 0$.

\medskip
More generally if $\omega \subset K$ is a Borel set let $\tau_{\omega}$ be the $n$-tuple of multiplication operators by the coordinate functions in $L^2(\omega,H_p \mid \omega)$ (this is the same as $\tau_K \mid L^2(\omega,H_p \mid \omega)$ since $L^2(\omega,H_p \mid \omega) \subset L^2(K,H_p \mid K))$. A key part of the proof of the main theorem will be to evaluate $k^-_p(\tau_{\omega})$ for increasingly general $\omega$, along lines similar of Lebesgue measure on ${\mathbb R}^n$ considered in \cite{7}.

We also define a constant $\gamma_K = \frac {(k^-_p(\tau_K))^p}{H_p(K)}$ where $p$ is the Hausdorff dimension of $K$. Lemma~$5.1$ and Lemma~$5.3$ imply that $0 < k^-_p(\tau_K) < \infty$ so that $0 < \gamma_K < \infty$.

\medskip
\noindent
{\bf Theorem 5.1.} 
{\em Let $\tau$ be a $n$-tuple of commuting Hermitian operators with $\sigma(\tau) \subset K$ and assume $p > 1$. Then we have
\[
(k^-_p(\tau))^p = \gamma_K \int_K m(x)dH_p(x)
\]
where $m$ is the multiplicity function of $\tau$.
}

\begin{proof}
Using Lemma~$5.2$ and the decomposition $\sH = \sH_{psing} \oplus \sH_{pac}$ the proof reduces to the case when the spectral measure of $\tau$ is absolutely continuous w.r.t.\ $H_p$, that is when $\sH = \sH_{pac}$. In view of Proposition~$2.1$ a further reduction is possible to the case when $\tau$ has finite cyclicity, that is when the multiplicity function is bounded. Since when $\tau$ has a cyclic vector and $H_p$-absolutely continuous spectral measure it is unitarily equivalent to a $\tau_{\omega}$, it means that the proof reduces to the case when $\tau = \tau_{\omega_1} \oplus \dots \oplus \tau_{\omega_m}$ for some Borel sets $\omega_j \subset K$, $1 \le j \le m$. In view of the last assertion in Proposition~$2.1$ the theorem holds for $\tau_{\omega_1} \oplus \dots \oplus \tau_{\omega_m}$ iff it holds for
\[
\tau_{F_{w_1}(\omega_1)} \oplus \dots \oplus \tau_{F_{w_m}(\omega_m)} 
\]
where $|w_1| = \dots = |w_m|$ because
\[
\tau_{F_{w_j}(\omega_j)} \simeq F_{w_j}(\tau_{\omega_j}).
\]
We may then choose $|w_j|$ sufficiently large and so that the $F_{w_j}(\omega_j)$, $1 \le j \le m$ are disjoint, which implies that
\[
\tau_{F_{w_1}(\omega_1)} \oplus \dots \oplus \tau_{F_{w_m}(\omega_m)} \simeq \tau_{\omega}
\]
where
\[
\omega = F_{w_1}(\omega_1) \cup \dots \cup F_{w_m(\omega_m)}.
\]
Thus the proof has been reduced to showing that
\[
(k^-_p(\tau_{\omega}))^p = \gamma_KH_p(\omega).
\]

First, assume $\omega$ is a finite union of $K_w$. Since $K_w$ is a disjoint union of $K_{w'}$, with $|w'| \ge |w|$ we may assume
\[
\omega = K_{w_1} \cup \dots \cup K_{w_m}
\]
where $|w_1| = \dots = |w_m|$ and $w_1,\dots,w_m$ are distinct. These $K_{w_j}$ are congruent and using again the last assertion in Prop.~$2.1$ the proof of this case reduces to proving the theorem for $\tau = \tau_{K_w} \otimes I_m$. The multiplicity function is $m$ times the indicator function of $K_w$ so that the right-hand side in the formula we want to prove is
\[
\begin{aligned}
&\gamma_KmH_p(K_w) \\
&= \frac {(k^-_p(\tau_K))^p}{H_p(K)} \cdot m \cdot \lambda^{p|w|} \cdot H_p(K) \\
&= m \cdot (\lambda^{|w|}k_p^-(\tau_K))^p = m(k^-_p(\tau_{K_w}))^p.
\end{aligned}
\]
On the left-hand side we have
\[
(k^-_p(\tau_{K_w} \otimes I_m))^p = (m^{1/p}k^-_p(\tau_{K_w}))^p
\]
by Thm.~$3.1$ which equals the right-hand side.

Next we prove the theorem for $\tau_{\omega}$ when $\omega \subset K$ is a general open subset. Let $\omega^{(L)}$ be the union of the $K_w \subset \omega$ with $|w| \le L$. The $\omega^{(L)}$ are clopen subsets of $K$ and are finite unions of $K_w$ so that the theorem holds for $\tau_{\omega^{(L)}}$ and the theorem for $\tau_{\omega}$ is obtained using Prop.~$2.1$, which gives $k^-_p(\tau_{\omega^{(L)}}) \uparrow k^-_p(\tau_{\omega})$.

Finally let $\omega \subset K$ be a Borel set and let $C$ be compact and $G$ in $K$ be open so that $C \subset \omega \subset G$ and $H_p(G\backslash C) < \epsilon$ for a given $\epsilon > 0$. We have $|k^-_p(\tau_{\omega}) - k^-_p(\tau_G)| \le |k^-_p(\tau_{G\backslash C})| = (\gamma_k\epsilon)^{1/p}$ using the fact that $G\backslash C$ is open in $K$ and Prop.~$2.1$. Thus
\[
\begin{aligned}
&|k^-_p(\tau_{\omega}) - (\gamma_KH_p(\omega))^{1/p}| \\
&\le |k^-_p(\tau_{\omega}) - k^-_p(\tau_G)| + |k^-_p(\tau_G) - (\gamma_KH_p(\omega))^{1/\omega}| \\
&\le (\gamma_K\epsilon)^{1/p} + |(\gamma_KH_p(G))^{1/p} - (\gamma_KH_p(\omega))^{1/p}| \\
&\le (\gamma_K\epsilon)^{1/p} + |(\gamma_K(H_p(\omega)+\epsilon))^{1/p} - (\gamma_KH_p(\omega))^{1/p}|.
\end{aligned}
\]
Since $\epsilon > 0$ was arbitrary, we get $k^-_p(\tau_{\omega}) = (\gamma_KH_p(\omega))^{1/p}$.
\end{proof}

\medskip
\noindent
{\bf Corollary 5.1.}
{\em Assume $\sigma(\tau) \subset K$ and $p > 1$. Then $k^-_p(\tau) = 0$ iff the spectral measure of $\tau$ is singular w.r.t.\ $H_p$.
}

\medskip
\noindent
{\bf Remark 5.1.}
In \cite{8} we showed for a $n$-tuple $\tau$ and a normed ideal $\sJ$ that there is a largest reducing subspace for $\tau$ on which $k_{\sJ}$ vanishes. In the case of commuting $n$-tuples of Hermitian operators and $\sJ = \sC^-_n$ this subspace is the subspace where the spectral measure is singular w.r.t.\ Lebesgue measure. The theorem we proved in this section shows that if $\sigma(\tau) \subset K$ and $p > 1$, then the largest reducing subspace on which $k^-_p$ vanishes for the restriction of $\tau$ is precisely $\sH_{psing}$.

\section{Concluding remarks}
\label{sec6}

\medskip
\noindent
{\bf Remark 6.1.}
It is natural to wonder whether in general
\[
(k^-_p(\tau_1 \oplus \tau_2))^p = (k^-_p(\tau_1))^p + (k_p^-(\tau_2))^p
\]
which would be much more than the ampliation homogeneity we proved. If $p = 1$ this is known to be true \cite{7}. For $1 < p \le \infty$ this is an open problem. While a negative answer would not be surprising, it is certainly desirable to clarify this issue.

\medskip
\noindent
{\bf Remark 6.2.}
To get results for more general self-similar fractals than the Cantor-like $K$ we considered it may be useful to replace $k^-_p$ by ${\tilde k}_p^-$ the variant of $k_p^-$ considered in \cite{7} pages 13--16. This amounts to extending the norms of normed ideals to $n$-tuples, not by the max of norm on the components but by the norm of $(T^*_1T_1 + \dots + T_n^*T_n)^{1/2}$ that is the modulus in the polar decomposition of the column $\left( \begin{matrix}
T_1 \\
\vdots \\
T_m
\end{matrix}
\right)$. This $|\tau|_{\sJ}^{\sim}$ has the advantage over $|\tau|_{\sJ}$ of being invariant under rotations, that is if $(u_{ij})_{1 \le i,j \le n}$ is a unitary matrix en the $n$-tuple $\left(\sum_j u_{ij}T_j\right)_{1 \le i \le n}$ has the same $\sim\sJ$-norm as $\tau = (T_i)_{1 \le i \le n}$. In particular ${\tilde k}^-_p$ may be better suited to handle self-similar sets $K$ when we use more general $F_i(x) = \lambda U_i(X - b(i)) + b(i)$ where $U_i \in O(n)$. In particular it is quite straightforward to use $\sim$-norms in \S 3 and to see that ampliation homogeneity still holds for ${\tilde k}_p^-$ which we record as the next theorem.

\medskip
\noindent
{\bf Theorem 6.3.}
{\em If $\tau$ is a $n$-tuple of bounded operators and $1 \le p \le \infty$ then
\[
{\tilde k}^-_p(\tau \otimes I_m) = m^{1/p}{\tilde k}_p^-(\tau).
\]
}

\medskip
\noindent
{\bf Remark 6.4.}
In \cite{10} we give an extension in another direction to the formula for $k^-_n(\tau)$ in \cite{7} to hybrid perturbations.


\end{document}